\documentclass[10pt]{article}   
\usepackage{epsfig,graphics}
\usepackage{cite}
\usepackage[T1]{fontenc}
\usepackage[swedish,english]{babel}
\usepackage{ntheorem}

\newcommand{\R}{{\bf R}}

\newcommand{\card}{\rm card}

\newtheorem{theorem}{Theorem}[section]

\newtheorem{lemma}{Lemma}[section]

\begin{document}

\title{\textsf{Almost sure convergence of the minimum bipartite matching
functional in Euclidean space}}

\author{\textsf{J.H.~Boutet de Monvel$^*$ and O.C.~Martin$^\dag$} \\ 
\textsf{\small $^*$Center for Hearing and Communication Research, Karolinska Institutet, 17176}\\ 
\textsf{\small Stockholm, Sweden; $^\dag$Laboratoire de Physique Th\'eorique et Mod\`eles Statistiques,}\\
 \textsf{\small Universit\'e de Paris-Sud, 91405 Orsay, France;}}

\date{To appear in Combinatorica}
\maketitle

\begin{abstract}
Let $L_N = L_{MBM}(X_1,\ldots ,X_N; Y_1,\ldots ,Y_N)$ be the minimum length of a
bipartite matching between two sets of points in $\mathbf{R}^d$, where
$X_1,\ldots ,X_N,\ldots$ and $Y_1,\ldots ,Y_N,\ldots$ are random points independently and 
uniformly distributed in $[0,1]^d$. We prove that for $d \ge 3$,  $L_N/N^{1-1/d}$ converges 
with probability one to a constant $\beta_{MBM}(d)>0$ as $N\to \infty $. 
\end{abstract}

\section{Introduction and statement of the result.}

\noindent Given two sets of $N$ points $X=\{X_1,...,X_N\}$ and $Y=\{Y_1,...,Y_N\}$ in
$\R^d$, a bipartite matching of $X$ and $Y$ is a perfect matching $M$ on the set $X\cup Y$, 
such that each pair in $M$ is made of one point of $X$ and one point of $Y$. The length of such 
a matching is defined to be the sum of the euclidean lengths of the edges formed by its pairs.
The (euclidean) minimum bipartite matching problem (MBMP) then asks one to find a 
bipartite matching of $X$ and $Y$ whose length is as small as possible. We shall denote by 
$L_{MBM}(X,Y)$ the length of a minimum bipartite matching of $X$ and $Y$.

A related problem is the simple minimum matching problem (MMP), where one is asked 
to find a perfect matching of smallest euclidean length on a set $X=\{X_1,...,X_N\}\subset \R^d$.
The subadditive methods inaugurated by Beardwood, Halton and Hammersley 
(BHH) \cite{BHH59_PCPS} and further developed 
in \cite{Steele81_AP,Rhee93_AAP,RedmondYukich94_AAP}, show 
that a strong limit theorem applies to the length $L_{MM}(X)$ of a simple minimum matching 
on $X$, when the points $X_1,\ldots, X_N$ are random.  
The theorem states that for any dimension $d$, if  $X_1,\ldots, X_N,\ldots$ is a sequence of 
points distributed independently and uniformly in a bounded region $\Omega\subset {\mathbf R}^d$, 
then the ratio $L_{MM}(X_1,\ldots X_N)/N^{1-1/d}$ converges almost surely to 
${\rm Vol(\Omega)}^{1/d}\beta_{MM}(d)$, where ${\rm Vol(\Omega)}$ denotes  the Lebesgues 
measure of $\Omega$ and $\beta_{MM}(d)>0$ is a universal constant depending only upon $d$.

The functional $L_{MBM}$ does not satisfy this form of limit theorem in dimensions 
$1$ and $2$. For $d=1$, the MBMP amounts to a sorting problem and it is not difficult 
to show that if $X$ and $Y$ both consist of $N$ points independently and uniformly 
distributed in $[0,1]$, there are constants $0<C_1<C_2$ such that 
$C_1\sqrt N\le L_{MBM}(X,Y)\le C_2 \sqrt N$ with probability $1-o(1)$ as
$N\to \infty$. Moreover in that case the variance of $L_{MBM}(X,Y)/\sqrt{N}$ does 
{\it not} converge to zero as $N\to \infty$. ($L_{MBM}$ is not ``self-averaging'', 
in the statistical physics' terminology.)
For $d=2$ Ajtai et al. \cite{Ajtai&Al84_C} proved a remarkable fact: if the sets 
$X,Y$ are now distributed in $[0,1]^2$, then for some constants $C_1,C_2$ indendent of 
$N$, one has $C_1\sqrt{N\log N}\le L_{MBM}(X,Y)\le C_2\sqrt{N\log N}$ with 
probability $1-o(1)$. Numerical simulations suggest that $L_{MBM}(X,Y)/\sqrt{N\log N}$ 
converges to a non-random constant as $N\to \infty$, however this has not yet been proved.

In this article, we show that for any $d\ge 3$ we recover a BHH theorem for the functional 
$L_{MBM}$.

\begin{theorem}\label{th1}
Let $X_1,...,X_N,...$ and $Y_1,...,Y_N,...$ be two sequences of
random points independently and uniformly distributed in $[0,1]^d$, where 
$d\ge 3$, and let $L_N = L_{MBM}(X_1,\ldots ,X_N;Y_1,\ldots ,Y_N)$.
There exists a constant $\beta_{MBM}(d)>0$ such that 
with probability one
$$ \lim_{N\to \infty} L_N/ N^{1-1/d} = \beta_{MBM}(d).$$
\end{theorem}

\section{Proof of Theorem \ref{th1}.}

To begin, we remark that to prove this theorem it will suffice to
establish that $L_N/N^{1-1/d}$ converges in mean value to a constant
$\beta_{MBM}(d)$. This is a consequence of the following lemma  \cite{Talagrand92_AAP}:

\begin{lemma}
For any $t>0$, one has
$$P(|{L_N\over N^{1-1/d}}- E({L_N\over N^{1-1/d}})| > t) \le 2 \exp(-{N^{1-2/d} t^2\over 8d}).$$
\end{lemma}

\noindent This result follows from the application of Azuma's inequality \cite{Azuma67_TMJ}  
and the martingale difference method to $L_N$, in a way by now standard in the
probabilistic theory of combinatorial optimisation \cite{Steele97_Book}.
Given the lemma, the theorem follows easily from the convergence of
$EL_N/N^{1-1/d}$ as $N\to \infty$, by applying the Borel-Cantelli lemma.

We have now to establish that for $d\ge 3$ the quantity
$EL_N/N^{1-1/d}$ indeed converges to a constant $\beta_{MBM}(d)>0$.
To prove this we exploit the subadditivity properties of $L_{MBM}$, in the spirit 
of Steele's theory of subadditive Euclidean functionals \cite{Steele81_AP}. 
Let us divide the unit cube $[0,1]^d$ into disjoint
similar subcubes $Q_k,~k=1,\ldots ,m^d$ with edges of length $1/m$,
and compare the value of $L_{MBM}(X,Y)$ to 
the sum 
\begin{equation} \label{SumOnCubes}
\sum_{k=1}^{m^d} L_k,
\end{equation}
where $L_k$ is the value of the functional $L_{MBM}$ for the set of points
$X_i$ and $Y_i$ which belongs to $Q_k$. A difficulty arises as in
general the $Q_k$'s do not contain the same number of points $X_i$ and of
points $Y_i$. (In fact the special properties of the MBMP in dimensions $1$ and 
$2$ originate from the fluctuations of the differences between these numbers 
around their mean value $0$.)
To give meaning to the sum (\ref{SumOnCubes}) we need to generalize the 
functional $L_{MBM}$ to matchings between two sets of different cardinalities. 
There are several ways to do this; we shall define 
$L_{MBM}(X_1,\ldots X_{N_1};Y_1,\ldots Y_{N_2})$  by imposing that 
the minimum matching contains as few unmatched points as possible.  That is if 
$N_1>N_2$, we leave $N_1-N_2$ points of $X$ unmatched, whereas if
$N_1<N_2$ we leave $N_2-N_1$ points of $Y$ unmatched. 

Although expression (\ref{SumOnCubes}) now makes sense, it is still not possible 
to write a subadditivity inequality of the same form as the one studied 
in \cite{Steele81_AP}. Indeed, such a form (which Steele calls ``geometric 
subadditivity'') implies an upper bound of the form $CN^{1-1/d}$ for the functional 
at hand \cite{Steele97_Book}, and it is easy to see that no such bound applies 
to $L_{MBM}(X,Y)$.  We shall however see that a geometric subadditivity 
property holds {\it in the mean} for the functional $L_{MBM}$.
Suppose that the points $X_1,\ldots X_{N_1},Y_1,\ldots Y_{N_2}$ belong to an
arbitrary cube $Q$ having edge length $a$, and divide $Q$ into
disjoint cubes $Q_p,~p=1,\ldots 2^d$ by splitting each edge in two halves.
Construct in each $Q_p$ an optimal matching in the sense just defined,
between the $n_{1,p}$ points $X_i$ and the $n_{2,p}$ points $Y_i$ in $Q_p$,
and denote its length by $L_p$.
The points that are left unpaired are in number $|n_{1,p}-n_{2,p}|$ in each
$Q_p$, so if $L_0$ denotes the length of an optimal matching for these
points one has
\begin{eqnarray} \label{Decimation}
L_{MBM}(X_1,\ldots X_{N_1};Y_1,\ldots ,Y_{N_2}) \le 
\sum_{p=1}^{2^d} L_p + L_0 \nonumber\\
\le \sum_{p=1}^{2^d} L_p + {1\over 2} a\sqrt d 
\sum_{p=1}^{2^d} |n_{1,p}-n_{2,p}|,
\end{eqnarray}
where the last inequality is obtained by bounding $L_0$ in an obvious way.

We shall apply this to $Q=[0,1]^d$. Let  $Q_{p_1}~p_1=1,\ldots 2^d$
be the cubes obtained in the above subdivision; let $Q_{p_1p_2}$ be
the cubes obtained by splitting in two halves the edges of each cube $Q_{p_1}$;
and so on. By repeating this operation $K$ times, we get a subdivision with
cubes $Q_{p_1\ldots p_K}$ whose edges are of length $1/2^K$. Let
$n_{1,p_1\ldots p_K}$ and $n_{2,p_1\ldots p_K}$ be respectively
the number of points $X_i$ and $Y_i$ in $Q_{p_1\ldots p_K}$. Apply
(\ref{Decimation}) first to the $Q_{p_1,\ldots p_{K-1}}$'s, then to the
$Q_{p_1\ldots p_{K-2}}$'s, etc, keeping at each step only those points which
are still unpaired. It is easy to convince oneself that the number of unpaired
points in each $Q_{p_1,\ldots p_{K-k}}$ just after step $k$ is given by
$|n_{1,p_1,\ldots p_{K-k}}-n_{2,p_1,\ldots p_{K-k}}|$. After step $k=K$ one
obtains a matching between $X_1,\ldots X_{N_1}$ and $Y_1,\ldots Y_{N_2}$
where all the points but $|N_1-N_2|$ are matched.
One is thus led to the following inequality:
\begin{eqnarray} \label{SousAddMBMP}
L_{MBM}(X_1,\ldots X_{N_1};Y_1,\ldots Y_{N_2}) 
\le \sum_{p_1\ldots p_K} L_{p_1\ldots p_K} \nonumber\\ 
+ \sum_{k=1}^K {\sqrt d\over 2^k}
\sum_{p_1\ldots p_k} |n_{1,p_1\ldots p_k}-n_{2,p_1\ldots p_k}|.
\end{eqnarray}
We now proceed to derive a subadditivity property for the mean 
value of $L_{MBM}(X,Y)$. We first consider the case where 
$N_1=\card X$ and $N_2=\card Y$ are not fixed integers but are independent Poisson 
random  variables with the same mean value $N$, the elements of $X$ and $Y$ being 
chosen independently and uniformly in $[0,1]^d$. For a given $k$, the numbers
$n_{1,p_1,\ldots p_k}$ and $n_{2,p_1,\ldots p_k}$ are then also independent
Poisson random variables, with parameter $N/2^{kd}$. Let 
$M(N)= EL_{MBM}(X_1,\ldots X_{N_1};Y_1,\ldots Y_{N_2})$. 
It is immediate by homogeneity that  
\begin{equation}
EL_{p_1\ldots p_K} = 2^{-K} M(N/2^{Kd}).
\end{equation}
Moreover from the well known properties of Poisson variables we have
\begin{equation} \label{RMSPoisson}
E|n_{1,p_1\ldots p_k}-n_{2,p_1\ldots p_k}| \le 
\sqrt 2 \Big( {N\over 2^{kd}} \Big)^{1/2}.
\end{equation}
By taking mean values in (\ref{SousAddMBMP}) we obtain:
\begin{equation}
M(N) \le 2^{K(d-1)}M(N/2^{Kd}) + \sqrt{2dN} \sum_{k=1}^K 2^{k(d/2-1)}.
\end{equation}
This inequality has been obtained for a subdivision of $[0,1]^d$ which
consists in $2^{Kd}$ similar cubes. Suppose now that we start from the
subdivision $\Sigma$ in $m^d$ similar cubes $Q_k~k=1,\ldots m^d$,
where $m$ is an arbitrary integer. One can then reproduce the previous
construction in the following manner. Let $m=2^K+r$
where $0\le r<2^K$. Consider the cube $Q_0=[0,2^{K+1}/m]^d$ and form the
natural subdivision $\Sigma_0$ of $Q_0$ by $2^{(K+1)d}$ cubes
$Q_{p_0,\ldots p_K}$ whose edges have length $1/m$. We can proceed with
$Q_0$ and $\Sigma_0$ to a $K+1$ steps construction similar to the one
which led to (\ref{SousAddMBMP}). The only differences are that $Q_0$ has
edges of length $2^{K+1}/m$ rather than $1$, and that some of the
$Q_{p_0\ldots p_K}$'s, namely those which belong to
$\Sigma_0$ but not to $\Sigma$, are empty. 
Nevertheless, we may write
\begin{eqnarray}
L_{MBM}(X_1,\ldots X_{N_1};Y_1,\ldots ,Y_{N_2}) - \sum_{p=1}^{m^d} L_k \nonumber \\ 
\le \sum_{k=0}^K {\sqrt d 2^{K-k} \over m}
\sum_{p_0\ldots p_k} |n_{1,p_0\ldots p_k}-n_{2,p_0\ldots p_k}| \nonumber \\
\le \sum_{k=0}^K {\sqrt d \over 2^k}
\sum_{p_0\ldots p_k} |n_{1,p_0\ldots p_k}-n_{2,p_0\ldots p_k}|.
\end{eqnarray}
Now $n_{1,p_0\ldots p_k}$ and $n_{2,p_0\ldots p_k}$ are Poisson 
variables with parameter lower than $2^{(K-k)d} N/m^d \le 2^{-kd}N$ so we 
still have
\begin{equation}
E|n_{1,p_0\ldots p_k}-n_{2,p_0\ldots p_k}| \le 
\sqrt 2 \Big({N \over 2^{kd}} \Big)^{1/2}.
\end{equation}
Taking average values one is led to 
\begin{equation}
M(N) \le m^{d-1} M(N/m^d) + 
2^d \sqrt{2dN} \sum_{k=0}^K 2^{k(d/2-1)}.
\end{equation}
Dividing this last inequality by $N^{1-1/d}$ and then replacing $N$ by
$m^dN$, we get
\begin{equation} 
{M(m^dN) \over (m^dN)^{1-1/d}} \le {M(N)\over N^{1-1/d}} + 
{2^d\sqrt{2d} \over N^{1/2-1/d}} \sum_{k=0}^K 2^{-k(d/2-1)}.
\end{equation}
If $d>2$, the sum on the r.h.s. of the last inequality is bounded above independently of 
$N$, and is divided by a positive power of $N$. Elementary analysis now shows that the 
ratio $M(N)/N^{1-1/d}$ necessarily converges to a limit $\beta_{MBM}(d)$ as $N\to \infty$. 
Indeed, let $f(t) = M(t^d)/t^{d-1}$. One verifies at once that $f(t)$ satisfies 
\begin{equation} \label{fInequality}
f(mt)\le f(t)+C/t^{d/2-1}
\end{equation}
for all $t>0$ and any integer $m$; $f(t)$ is continuous,
since $M(N)$ is a continuous function of $N$.
So the expression $f(t) + C_d/t^{d/2-1}$ is bounded in $[1,2]$ and since
$[1,\infty[$ is the union of the intervals $m[1,2], m\ge 1$, it follows
from (\ref{fInequality}) that $f(t)$ remains bounded as $t\to \infty$,
thus $\lim^* f(t) < \infty$. Now define $\beta=\lim_* f(t)$. For any
$\epsilon >0$, chose $t_0\gg 1$ and
$\eta >0$ such that $f(t)+C_d/t^{d/2-1} < \beta + \epsilon$
for $t$ in the interval $I=[t_0-\eta,t_0+\eta]$.
Since the intervals $mI$, $m\ge 1$ span a whole interval
$[A,\infty[$ for an $A$ sufficiently large,
it follows again from (\ref{fInequality}) that
$\lim^* f(t)\le \beta+\epsilon$.
Since $\epsilon$ is arbitrary one has $\lim^* f(t)=\beta$, hence
$f(t) \to \beta$ as $t\to \infty$,  from which it follows that 
$\lim_{N\to \infty} M(N)/N^{1-1/d}=\beta$. Q.E.D.

We have thus shown for $d\ge 3$, that one has
\begin{equation} \label{PoissonMBMPAsymptotics}
EL_{MBM}(X_1,\ldots,X_{N_1};Y_1,\ldots,Y_{N_2}) 
\sim \beta_{MBMP}^E(d)N^{1-1/d},~N\to \infty
\end{equation}
when $N_1$ and $N_2$ are independent Poisson variables with parameter $N$.
The same result for the mean value $EL_N$, where $N$ is a fixed integer,
follows then easily. Indeed, we have the obvious bound
\begin{eqnarray}
|L_{MBM}(X_1,\ldots X_N;Y_1,\ldots Y_N)-
L_{MBM}(X_1,\ldots X_{N_1};Y_1,\ldots Y_{N_2})| \nonumber\\
\le \sqrt d (|N_1-N|+|N_2-N|),
\end{eqnarray}
whence taking mean values, 
\begin{equation}
|EL_N - EL_{MBM}(X_1,\ldots X_{N_1};Y_1,\ldots Y_{N_2})| \le 2 \sqrt{2dN},
\end{equation}
and dividing by $N^{1-1/d}$ we deduce that 
\begin{equation}
\lim_{N\to \infty} {EL_N\over N^{1-1/d}} \to \beta_{MBM}(d).
\end{equation}
Theorem \ref{th1} is now proved.

\section{Concluding remarks.}

\noindent 1) Our decimation procedure does not give back the bounds
proven by  Ajtai {\it et al.} in $d=2$, but a weaker 
$O(\sqrt{N} \ln N)$ bound. 
It is believed that a self-averaging theorem applies also to the 
functional $L_{MBM}$ in dimension $2$ \cite{Smith89_Thesis}.

\noindent 2) The estimation of the constants $\beta_{MBM}(d)$ is also an
interesting problem. A remarkable result of Talagrand \cite{Talagrand92_AAP}
shows that one has $\beta_{MBM}(d)= \sqrt{d/2e\pi} (1+O(\ln d / d))$ as 
$d\to \infty$. It is conjectured that a $1/d$ series expansion actually exists 
for $\beta_{MBM}(d)$.

\noindent 3) M\'ezard and Parisi have obtained detailed analytic predictions for  
the {\it random link} versions of the MMP and the MBMP  \cite{MezardParisi87_JdP}, 
where the distance matrix between the points $X_i$ and $Y_j$ is replaced by a matrix of 
independent and identically distributed entries. (Some of these predictions, for the random 
assignment problem, have been proven recently by Aldous \cite{Aldous01_RSA}.) 
Numerical studies \cite{BoutetMartin97_PRL,HBM98_EPJB} indicate that for the MMP and the 
MBMP, the random link model provides one with a very good ``mean-field'' approximation to 
the Euclidean model in the large $d$ limit. Except for simpler combinatorial problems 
however \cite{BertsimasVanRyzin90_ORL}, very few rigorous results are known for comparing 
the euclidean and the random link models.

\bigskip
{\noindent \bf \large Aknowledgments}

\noindent It is a pleasure to thank J.M. Steele for fruitful discussions and pointing to us 
reference \cite{Talagrand92_AAP}.

\bibliography{co,jbdm}
\bibliographystyle{perroten}

\end{document}